\documentclass[11pt, reqno]{amsart}

\usepackage{amsthm,amsmath,amssymb,xcolor,tikz}

\usetikzlibrary{decorations.pathreplacing}
\usetikzlibrary{decorations.markings}
\usetikzlibrary{arrows}

\usepackage[margin=1in]{geometry}
\usepackage{mathtools}

\usepackage[pdfstartview=FitV,hidelinks,colorlinks=true]{hyperref}

\usepackage[utf8]{inputenc}
\usepackage[T1]{fontenc}

\usepackage[hidelinks]{hyperref}
\usepackage{url}

\usepackage[noabbrev,capitalize]{cleveref}

\usepackage[color,final]{showkeys} %add in 'final' into parameter to remove showkeys

\colorlet{refkey}{orange!20}
\colorlet{labelkey}{blue!30}
\newtheorem{theorem}{Theorem}

\newtheorem*{question*}{Question}

\theoremstyle{definition}

\newtheorem*{definition*}{Definition}

\theoremstyle{remark}

\newcommand{\paren}[1]{\left(#1\right)}

\newcommand{\abs}[1]{\left\lvert#1\right\rvert}

\newcommand{\eps}{\varepsilon}

\newcommand{\mb}{\mathbb}

\title{Paths of given length in tournaments}

\author[Sah]{Ashwin Sah}
\author[Sawhney]{Mehtaab Sawhney}
\author[Zhao]{Yufei Zhao}
\address{Massachusetts Institute of Technology, Cambridge, MA 02139, USA}
\email{\{asah,msawhney,yufeiz\}@mit.edu}
\date{}

\thanks{
Sah and Sawhney were supported by NSF Graduate Research Fellowship Program DGE-1745302. Zhao was supported by NSF Award DMS-1764176, the MIT Solomon Buchsbaum Fund, and a Sloan Research Fellowship.
}

\begin{document}

\begin{abstract}
We prove that every $n$-vertex tournament has at most $n\paren{\frac{n-1}{2}}^k$ walks of length $k$.
\end{abstract}

\maketitle
We determine the maximum density of directed $k$-edge paths in an $n$-vertex tournament.
Our focus is on the case of fixed $k$ and large $n$.
The expected number of directed $k$-edge paths in a uniform random $n$-vertex tournament is $n(n-1)\cdots(n-k)/2^k = (1+o(1))n(n/2)^k$.
In this short note we show that one cannot do better, thereby confirming an unpublished conjecture of Jacob Fox, Hao Huang, and Choongbum Lee.
The \emph{length} of a path or walk refers to its number of edges.

\begin{theorem} \label{thm:main}
Every $n$-vertex tournament has at most $n\paren{\frac{n-1}{2}}^k$ walks of length $k$.
\end{theorem}

Every regular tournament (with odd $n$) has exactly $n\paren{\frac{n-1}{2}}^k$ walks of length $k$, thereby attaining the upper bound in the theorem. 
On the other hand, the transitive tournament minimizes the number of $k$-edge paths (or walks) among $n$-vertex tournaments. Indeed, a folklore result (with an easy induction proof) says that every tournament contains a directed Hamilton path. So every $(k+1)$-vertex subset contains a path of length $k$. Hence every $n$-vertex tournament contains at least $\binom{n}{k+1}$ paths of length $k$, with equality for a transitive tournament.

\medskip

\tikzstyle{v}=[draw, circle, fill, inner sep = 1pt]
\tikzstyle{e}=[postaction={decorate,decoration={
        markings,
        mark=at position .5*\pgfdecoratedpathlength+1.5pt 
        with {\arrow{angle 90}}
      }}]
\tikzstyle{br}=[decorate,decoration={brace,amplitude=5pt}]

Here is a  ``proof by picture'' of \cref{thm:main}. A more detailed proof is given later.
A different proof, using entropy, by
 Dingding Dong and Tomasz \'Slusarczyk, 
is given in the appendix. 
\begin{align*}
\bigl(
\begin{tikzpicture}[scale=.5, baseline=-.8ex]
	\draw[e] (0,0) node[v]{} -- ++(1,0) node[v]{};
\end{tikzpicture}
\bigr)^{-1}
\bigl(
\begin{tikzpicture}[scale=.5, baseline=-.8ex]
	\draw[e] (0,0) node[v]{} -- ++(1,0) node[v]{};
	\node[align=center] at (2,0) {$\cdots$};
	\path (3,0) node[v](v0){} 
		-- ++(1,0) node[v](v1){} 
		-- ++(1,0) node[v](v2){} 
		-- ++(1,0) node[v](v3){} 
		-- ++(1,0) node[v](v4){};
	\draw[e](v0)--(v1);
	\draw[e](v1)--(v2);
	\draw[e](v2)--(v3);
	\draw[e](v3)--(v4);
	\draw[br] (0,.3) -- ++(7,0) node[above=.5ex,midway] {\scriptsize $k$ edges};
\end{tikzpicture}
\bigr)^2
\
&\le
\
\begin{tikzpicture}[scale=.5, baseline=-.8ex]
	\draw[e] (0,.3) node[v]{} -- ++(1,0) node[v]{};
	\draw[e] (0,-.3) node[v]{} -- ++(1,0) node[v]{};
	\node[align=center] at (2,0) {$\cdots$};
	\path (3,.3) node[v](v0){} 
		-- ++(1,0) node[v](v1){} 
		-- ++(1,0) node[v](v2){} 
		-- ++(1,-.3) node[v](v3){} 
		-- ++(1,0) node[v](v4){};
	\path (3,-.3) node[v](w0){} 
		-- ++(1,0) node[v](w1){} 
		-- ++(1,0) node[v](w2){};
	\draw[e] (v0)--(v1);
	\draw[e] (v1)--(v2);
	\draw[e] (v2)--(v3);
	\draw[e] (v3)--(v4);
	\draw[e] (w0)--(w1);
	\draw[e] (w1)--(w2);
	\draw[e] (w2)--(v3);
	\draw[br](0,.6) -- ++(6,0) node[above=.5ex,midway] {\scriptsize $k-1$ edges};
\end{tikzpicture}
&&\text{\footnotesize Cauchy--Schwarz}
\\
\
&\le
\
\begin{tikzpicture}[scale=.5, baseline=-.8ex]
	\draw[e] (0,.6) node[v]{} -- ++(1,0) node[v]{};
	\draw[e] (0,0) node[v]{} -- ++(1,0) node[v]{};
	\node[align=center] at (2,0) {$\cdots$};
	\path (3,.6) node[v](v0){} 
		-- ++(1,0) node[v](v1){} 
		-- ++(1,-.3) node[v](v2){} 
		-- ++(1,-.3) node[v](v3){} 
		-- ++(1,0) node[v](v4){};
	\path (3,0) node[v](w0){} 
		-- ++(1,0) node[v](w1){};	
	\draw[e] (v0)--(v1);	
	\draw[e] (v1)--(v2);	
	\draw[e] (v2)--(v3);	
	\draw[e] (v3)--(v4);	
	\draw[e] (w0)--(w1);	
	\draw[e] (w1)--(v2);	
	\draw[e] (5,-.3) node[v](x2){} -- (v3);
	\draw[br] (0,.9) -- ++(5,0) node[above=.5ex,midway] {\scriptsize $k-2$ edges};
\end{tikzpicture}
&& \text{\footnotesize $2ab\le a^2 + b^2$}
\\
\ &\le \ 
\begin{tikzpicture}[scale=.5, baseline=-.8ex]
	\draw[e] (0,.3) node[v]{} -- ++(1,0) node[v]{};
	\draw[e] (0,-.3) node[v]{} -- ++(1,0) node[v]{};
	\node[align=center] at (2,0) {$\cdots$};
	\path (3,.3) node[v](v0){} 
		-- ++(1,0) node[v](v1){} 
		-- ++(1,-.3) node[v](v2){} 
		-- ++(1,0) node[v](v3){};
	\path (3,-.3) node[v](w0){} 
		-- ++(1,0) node[v](w1){};	
	\draw[e] (v0)--(v1);	
	\draw[e] (v1)--(v2);	
	\draw[e] (v2)--(v3);	
	\draw[e] (w0)--(w1);	
	\draw[e] (w1)--(v2);	
	\draw[br] (0,.6) -- ++(5,0) node[above=.5ex,midway] {\scriptsize $k-2$ edges};
\end{tikzpicture}
\ \cdot \paren{\tfrac{n-1}{2}}^2 \
&& \text{\footnotesize $d^+(x) d^-(x) \le \paren{\tfrac{n-1}{2}}^2$}
\\
\cdots&\textnormal{ Iterating}\cdots 
\\
\ & \le \ 
\begin{tikzpicture}[scale=.5, baseline=-.8ex]
	\draw[e] (0,0) node[v]{} -- ++(1,0) node[v]{};
\end{tikzpicture}
\ \cdot \paren{\tfrac{n-1}{2}}^{2k-2}.  \
\end{align*}

Let us mention some related problems and results.
The most famous open problem with this theme is Sidorenko's conjecture~\cite{ES83,Sid93}, which says that for a fixed bipartite graph $H$, among graphs of a given density, quasirandom graphs minimize $H$-density. 
For recent progress on Sidorenko's conjecture see~\cite{CKLL18,CL21}.

Zhao and Zhou~\cite{ZZ} determined all directed graphs that have constant density in all tournaments; they are all disjoint unions of trees that are each constructed in a recursive manner, as conjectured by Fox, Huang, and Lee. 
As discussed at the end of \cite{ZZ}, it would be interesting to characterize directed graphs $H$ where is the $H$-density in tournaments maximized by the quasirandom tournament (such $H$ is called \emph{negative}), and likewise when ``maximized'' is replaced by ``minimized'' (such $H$ is called \emph{positive}). 
Our main result here implies that all directed paths are negative. 
It would be interesting to see what happens for other edge-orientations of a path.
Starting with a negative (resp.~positive) digraph, one can apply the same recursive construction as in \cite{ZZ} to produce additional negative (resp.~positive) digraphs, namely by taking two disjoint copies of the digraph and adding a single edge joining a pair of twin vertices.

The problem of maximizing the number of directed $k$-cycles in a tournament is also interesting and not completely understood.
Recently, Grzesik, Kr\'al', Lov\'asz, and Volec~\cite{GKLV23} showed that quasirandom tournaments maximize the number of directed $k$-cycles whenever $k$ is not divisible by 4.
On the other hand, when $k$ is divisible by 4, quasirandom tournaments do not maximize the density of directed $k$-cycles.
The maximum directed $k$-cycle density is known for $k=4$~\cite{BH65,Col64} and $k=8$~\cite{GKLV23} but open for all larger multiples of $4$.
See \cite{GKLV23} for discussion.

A related problem is determining the maximum number $P(n)$ of Hamilton paths in a tournament (the problem for Hamilton cycles is related).
By considering the expected number of Hamilton paths in a random tournament, one has $P(n) \ge n!/2^{n-1}$.
This result, due to Szele~\cite{Sze43}, is considered the first application of the probabilistic method. 
This lower bound has been improved by a constant factor \cite{AAR01,Wor04}.
Alon~\cite{A90} proved a matching upper bound of the form $P(n) \le n^{O(1)} n!/2^{n-1}$ (also see \cite{FK05} for a later improvement).

%
%Now we prove \cref{thm:main}, which, by a standard graph limit argument, is equivalent to the following analytic statement. To deduce \cref{thm:main}, given an $n$ vertex tournament $T$ with vertex set $\{1,\ldots,n\}$, we apply \cref{thm:analytic} where the function $f(x,y)$ is $1$ when there is a directed edge from $n\lceil\frac{x}{n}\rceil$ to $n\lceil\frac{y}{n}\rceil$ in $W$, and $0$ otherwise.
%
%\begin{theorem}\label{thm:analytic}
%Suppose $f:[0,1]^2\to[0,1]$ is an integrable function satisfying $f(x,y) + f(y,x) \le 1$ for all $x,y \in [0,1]$. 
%Then for all positive integers $k$,
%\[\int_{[0,1]^{k+1}}f(x_0,x_1) f(x_1,x_2) \cdots f(x_{k-1},x_k) \, dx_0\cdots dx_k \le 2^{-k}.\]
%\end{theorem}

%\begin{remark}
%For $k\ge 2$, the proof shows that equality holds if and only if $\int_{[0,1]} f(x,y)dy= 1/2$ for almost all $x$, corresponding to limits of regular tournaments.
%\end{remark}

\begin{proof}[Proof of \cref{thm:main}]
We may assume that $k\ge 1$.
Let $f(x,y) = 1$ if $(x,y)$ is a directed edge in the tournament, and $0$ otherwise.
Let $g_t(x)$ denote the number walks of length $t$ ending at $x$. 
Let $d^+(x)$ and $d^-(x)$ denote the out-degree and the in-degree of $x$, respectively.
We have, for each $t \ge 1$,
\[
g_t(y) = \sum_x g_{t-1}(x) f(x,y).
\]
Define
\[
A_t 
\coloneqq 
\sum_y g_t(y)^2 d^+(y).
=
\sum_{y,z} g_t(y)^2 f(y,z).
\]
We have, for each $t \ge 1$,
\begin{align*}
A_t 
&= \sum_{x,x',y} g_{k-1}(x)f(x,y)g_{k-1}(x')f(x',y) d^+(y)
\\
&\le \sum_{x,x',y} \paren{\frac{g_{k-1}(x)^2 + g_{k-1}(x')^2}{2}}f(x,y)f(x',y) d^+(y)
\\
&= \sum_{x,x',y} g_{k-1}(x)^2 f(x,y)f(x',y) d^+(y)
\\
&= \sum_{x,y} g_{k-1}(x)^2 f(x,y)d^-(y)d^+(y)
%\\
%&\le \sum_{x,y} g_{k-1}(x)^2 f(x,y) \paren{\frac{d^-(y) + d^+(y)}{2}}^2
\\
&\le \paren{\frac{n-1}{2}}^2 \sum_{x,y} g_{k-1}(x)^2 f(x,y) \qquad \text{\small [since $\textstyle  d^-(y) d^+(y)  \le \paren{\frac{n-1}{2}}^2$]}
\\
&= \paren{\frac{n-1}{2}}^2 A_{t-1}.
\end{align*}
So, for all $t \ge 0$,
\[
A_t \le A_0\paren{\frac{n-1}{2}}^{2t} \le n \paren{\frac{n-1}{2}}^{2t+1}.
\]
Let $W_k$ be the number of walks of length $k$.
Applying the Cauchy--Schwarz inequality,
\begin{align*}
W_k = \sum_{y} g_{k-1}(y) d^+(y)
&
%\le \paren{\sum_y g_{k-1}(x)^2 d^+(y)}^{1/2} \paren{\sum_y d^+(y)}^{1/2}
\le \sqrt{\sum_y g_{k-1}(x)^2 d^+(y)} \sqrt{\sum_y d^+(y)}
\le \sqrt{A_{k-1} A_0}
\le n \paren{\frac{n-1}{2}}^k.
\qedhere
\end{align*}
\end{proof}

The above proof also gives the following stability result.
\begin{theorem}[Stability] \label{thm:stability}
For $k\ge 2$, an $n$-vertex tournament satisfying 
\[
\sum_x \abs{d^+(x) - \frac{n-1}{2}} 
\ge \eps \binom{n}{2}
\]
has at most $(1-\frac{\eps^2}{2})n\paren{\frac{n-1}{2}}^k$ walks of length $k$.
\end{theorem}
Note that by symmetry, we can replace $d^+$ by $d^-$ in the hypothesis of \cref{thm:stability}.
\begin{proof}
We use the notation from the earlier proof.
We have
\begin{align*}
W_2 = 
\sum_x d^+(x)d^-(x) 
&\le 
\sum_x d^+(x)(n-1 - d^+(x))
\\
&= 
\sum_x \paren{ \paren{\frac{n-1}{2}}^2 - \paren{\frac{n-1}{2} - d^+(x)}^2}
\\
&\le 
n\paren{\frac{n-1}{2}}^2 - \frac{1}{n} \paren{\sum_x \abs{\frac{n-1}{2} - d^+(x)}}^2
\\
&\le (1-\eps^2) n\paren{\frac{n-1}{2}}^2.
\end{align*}
From the proof of \cref{thm:main}, we have
\[
W_k^2 \le 
A_{k-1} A_0
\le \paren{\frac{n-1}{2}}^{2(k-2)} A_1 A_0.
\]
Using $A_0 \le n(n-1)/2$ and $A_1 = \sum_x d^-(x)^2 d^+(x)$, we obtain
\[
W_k^2
\le n \paren{\frac{n-1}{2}}^{2k-3} \sum_x d^-(x)^2 d^+(x).
\]
In the proof of \cref{thm:main}, we defined $g_k(x)$ to be the number of $k$-edge walks ending at $x$.
By running the same proof for the number of $k$-edge walks starting at $x$, we deduce
\[
W_k^2 
\le 
n \paren{\frac{n-1}{2}}^{2k-3} \sum_x d^-(x) d^+(x)^2.
\]
Taking the average of the two bounds, we obtain
\begin{align*}
W_k^2 
&\le
n \paren{\frac{n-1}{2}}^{2k-3}
\sum_x d^-(x) d^+(x) \paren{\frac{d^-(x) + d^+(x)}{2}}
\\
&\le 
n \paren{\frac{n-1}{2}}^{2k-2}
\sum_x d^-(x) d^+(x) 
\\
&\le (1-\varepsilon^2) n^2 \paren{\frac{n-1}{2}}^{2k}
\le \paren{\paren{1-\frac{\varepsilon^2}{2}} n \paren{\frac{n-1}{2}}^{k}}^2. \qedhere 
\end{align*}
\end{proof}

% \bibliographystyle{amsplain0}
% \bibliography{path-tournament}

\appendix
\section{An entropy proof 
\\ 
by~Dingding Dong and Tomasz \'Slusarczyk
}

Here is another proof of \cref{thm:main} using entropy.
Given a discrete random variable $X$ taking values in $\Omega$, its entropy is defined as
\[H(X) = -\sum_{x\in \Omega}\mb{P}(X=x)\log \mb{P}(X=x).\]
We have the uniform bound
\[H(X)\le\log|\Omega|.\]
The chain rule says that if $X$ and $Y$ are jointly distributed random variables, then
\[H(X,Y) = H(X) + H(Y|X),\]
where the conditional entropy $H(Y|X)$ is defined as
\[H(Y|X) = \sum_{x\in\Omega}\mb{P}(X=x)H(Y|X=x).\]
Here $H(Y|X=x)$ is the entropy of the conditional distribution of $Y$ given $X=x$.

\begin{proof}[Entropy proof of \cref{thm:main}]
Consider a random walk $X_1,\ldots,X_{k+1}$ chosen uniformly from the set of all $W_k$ walks of length $k$ in the given tournament.
This random walk is Markovian in the sense that the distribution of $(X_i,\ldots,X_{k+1})$ conditional on $(X_1,\ldots,X_i)$ is the same as the distribution of $(X_i,\ldots, X_{k+1})$ conditional on $X_i$. Indeed, this conditional distribution is uniform over all walks $(X_i,\ldots, X_{k+1})$ with a given starting vertex $X_i$. 
In particular, $H(X_j|X_{j-1},\ldots,X_1)=H(X_j|X_{j-1})$.

Applying the chain rule, we have 
\begin{align*}
\log W_k =
H(X_1,\ldots,X_{k+1}) &= H(X_1, X_2) + \sum_{j=2}^{k}H(X_{j+1}|X_1,\ldots,X_j)\\
&= H(X_1,X_2) + \sum_{j=2}^{k}H(X_{j+1}|X_{j}).
\end{align*}
Likewise,
\[H(X_1,\ldots,X_{k+1}) = H(X_{k+1},X_k) + \sum_{j=2}^{k}H(X_{j-1}|X_j).\]
Taking the average of the two bounds, we obtain
\[
H(X_1,\ldots,X_{k+1}) 
= \frac{H(X_1,X_2)+ H(X_{k+1},X_{k})}{2} + \frac{1}{2}\sum_{j=2}^k(H(X_{j-1}|X_j)+H(X_{j+1}|X_j)).
\]
For each $2 \le j \le k$ and vertex $x$,
by the uniform bound,
$H(X_{j-1}|X_j = x) \le \log d^-(x)$
and
$H(X_{j+1}|X_j = x) \le \log d^+(x)$.
Also, $d^-(x) d^+(x) \le (n-1)^2/4$. 
Thus
\begin{align*}
H(X_{j-1}|X_j = x) + H(X_{j+1}|X_j = x)
\le \log d^-(x) + \log d^+(x)
\le 2 \log  \paren{\frac{n-1}{2}}.
\end{align*}
Thus
\[
H(X_{j-1}|X_j)+H(X_{j+1}|X_j) \le 2 \log  \paren{\frac{n-1}{2}}.
\]
Also $H(X_j, X_{j+1}) \le \log \binom{n}{2}$ by the uniform bound.
Therefore
\[
\log W_k = H(X_1,\ldots,X_{k+1})
\le \log \binom{n}{2} + (k-1) \log \paren{\frac{n-1}{2}}=\log \paren{ n\paren{\frac{n-1}{2}}^k}. \qedhere
\]
\end{proof}

\end{document}